\documentclass[fontsize=11pt,bibliography=toc]{scrartcl}

\setkomafont{title}{\sffamily}
\setkomafont{pagehead}{\sffamily}
\setkomafont{author}{\sffamily}
\addtokomafont{author}{\vspace*{2em}\large}

\usepackage[automark]{scrlayer-scrpage}
\usepackage{mathtools}
\ihead{\headmark}
\ohead{\pagemark}
\usepackage{tikz}
\usetikzlibrary{arrows.meta,decorations.markings}
\definecolor{qblue}{RGB}{30,136,229}
\definecolor{cred}{RGB}{228,52,38}

\usepackage[
  paper=a4paper,
  left=80pt,
  right=80pt,
  top=86pt,
  bottom=56pt,
  bindingoffset=0mm,
  includefoot,
  includehead,
  foot=16.5pt,
  head=16.5pt
]{geometry}

\DeclareMathAlphabet{\pazocal}{OMS}{zplm}{m}{n}
\DeclareMathAlphabet{\mymathbb}{U}{bbold}{m}{n}

\usepackage{indentfirst}
\usepackage{lipsum}
\usepackage[english]{babel}
\usepackage{amsmath} % must load amsmath before amsthm
\usepackage{amsthm}
\usepackage{amssymb}
\usepackage[sfdefault]{FiraSans}
\usepackage[charter,expert]{mathdesign}
\usepackage[scaled=.96]{XCharter}% matches the size used in math
\usepackage{setspace}
\usepackage{nicefrac}% compact symbols for fractions
\usepackage[dvipsnames]{xcolor}
    \definecolor{linkcolor}{HTML}{0D6A9E}
    \definecolor{3blue}{HTML}{0072B2}
    \definecolor{3green}{HTML}{009E73}
    \definecolor{3ochre}{HTML}{E69F00}
    \definecolor{3yellow}{HTML}{F0E442}
    \definecolor{3cyan}{HTML}{56B4E9}
    \definecolor{3red}{HTML}{D55E00}
    \definecolor{3pink}{HTML}{CC79A7}
    \definecolor{2blue}{HTML}{1A85FF}
    \definecolor{2red}{HTML}{D41159}
    \usepackage{tikz}
\usetikzlibrary{arrows.meta,decorations.markings}
\definecolor{qblue}{RGB}{30,136,229}
\definecolor{cred}{RGB}{228,52,38}
\usepackage{graphicx}
    \graphicspath{ {./gfx/} }
\usepackage{caption}
\usepackage{tensor}
\usepackage{tikz-cd}
\tikzcdset{arrow style=tikz}

\usepackage{hyperref}
    \hypersetup{
    colorlinks=true,
    linktoc=all,
    linkcolor=linkcolor,
    citecolor=linkcolor,
    filecolor=magenta,      
    urlcolor=linkcolor,
    pdftitle={The monodromy manifold of the difference Painlev\'e equation of type $A_2^{(1)*}$},
    pdfpagemode=FullScreen,
    }
\usepackage{cleveref}

    \newtheorem{theorem}{Theorem}[section]

    \newtheorem{definition}[theorem]{Definition}
    \newtheorem{proposition}[theorem]{Proposition}

    \newtheorem{remark}[theorem]{Remark}
   \newtheorem*{conjecture}{Conjecture}

\newcommand{\diag}{\mathrm{diag}}

\begin{document}

\title{The monodromy manifold of the difference Painlev\'e equation of type $A_2^{(1)*}$ \vspace{-.5em}}
\newcommand{\orcidlink}[1]{%
  \raisebox{-.5ex}{\href{https://orcid.org/#1}{\includegraphics[height=15pt]{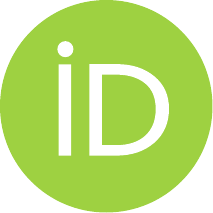}}}\hspace{.5em}}

\author{\orcidlink{0000-0001-9917-2547}Marta Mazzocco,
        \orcidlink{0009-0007-0946-7626}Anna Nefedova,
        \orcidlink{0000-0002-2765-0957}Volodya Rubtsov}
\date{\vspace{-2.5em}}

\maketitle
\begin{abstract}
   \noindent\textbf{Abstract.} We study the monodromy manifold of the difference Painlevé equation $d$-$P(A_2^{(1)*})$. Following Dzhamay, Sakai and Takenawa, this equation is realised as a composition of discrete Schlesinger transformations acting on $3\times3$ Fuchsian systems with three simple poles, so that, after a scalar normalisation, its monodromy data are points of the $\mathrm{SL}_3(\mathbb{C})$ character variety of the three-holed sphere. Starting from Lawton's presentation of this character variety, we show that, once the conjugacy classes of the local monodromies are fixed, it is given by a six-parameter family of affine cubic surfaces. One step of the discrete dynamics acts on this family as a central twist of order three. We parametrise the character variety by a triple cover of the Fock-Goncharov torus, in which the coordinates and the coefficients of the cubic are Laurent polynomials with positive coefficients. In a chewing-gum limit, in which the two Fock-Goncharov coordinates attached to one edge of the triangulation tend to infinity, the cubic degenerates to the monodromy manifold of the sixth Painlevé equation, and the traces of the four $\mathrm{SL}_2(\mathbb{C})$ local monodromies appear as explicit functions of the surviving coordinates. Based on these results, we conjecture explicit affine del Pezzo surfaces as monodromy manifolds of the difference Painlevé equations of types $A_1^{(1)*}$ and $A_0^{(1)**}$.
\end{abstract}
\section{Introduction}
\begin{figure}[h]
\centering
\begin{tikzpicture}[x=1cm,y=1cm,>=stealth,
  every node/.style={inner sep=2pt},
  eqname/.style={font=\scriptsize,inner sep=1pt}]

% ---------- elliptic ----------
\node[text=3green] (ell) at (0,4.4) {$A_0^{(1)}$};
\node[eqname,text=3green,anchor=west] at (ell.east) {elliptic};

% ---------- multiplicative (q-difference) ----------
\node[text=3blue] (m0)  at (0,2.8)    {$A_0^{(1)*}$};
\node[text=3blue] (m1)  at (1.6,2.8)  {$A_1^{(1)}$};
\node[text=3blue] (m2)  at (3.2,2.8)  {$A_2^{(1)}$};
\node[text=3blue] (m3)  at (4.8,2.8)  {$A_3^{(1)}$};
\node[text=3blue] (m4)  at (6.4,2.8)  {$A_4^{(1)}$};
\node[text=3blue] (m5)  at (8.2,2.8)  {$A_5^{(1)}$};
\node[text=3blue] (m6)  at (10.1,2.8) {$A_6^{(1)}$};
\node[text=3blue] (m7)  at (11.9,3.8) {$A_7^{(1)}$};
\node[text=3blue] (m7p) at (11.9,2.8) {$A_7^{(1)\prime}$};
\node[text=3blue] (m8)  at (13.6,2.8) {$A_8^{(1)}$};

\node[eqname,text=3blue,above] at (m3.north) {$q$-$PVI$};
\node[eqname,text=3blue,above] at (m4.north) {$q$-$PV$};
\node[eqname,text=3blue,above] at (m5.north) {$q$-$PIV$};
\node[eqname,text=3blue,above] at (m8.north) {$q$-$PI$};

% ---------- additive (difference) ----------
\node[text=3ochre] (a0) at (0,1.2)   {$A_0^{(1)**}$};
\node[text=3ochre] (a1) at (1.6,1.2) {$A_1^{(1)*}$};
\node[text=3ochre] (a2) at (3.2,1.2) {$A_2^{(1)*}$};

% ---------- differential ----------
\node[text=3red] (d4) at (4.8,1.2)   {$D_4^{(1)}$};
\node[text=3red] (d5) at (6.4,1.2)   {$D_5^{(1)}$};
\node[text=3red] (d6) at (8.2,1.2)   {$D_6^{(1)}$};
\node[text=3red] (d7) at (10.1,1.2)  {$D_7^{(1)}$};
\node[text=3red] (d8) at (11.9,1.2)  {$D_8^{(1)}$};
\node[text=3red] (e6) at (8.2,-0.4)  {$E_6^{(1)}$};
\node[text=3red] (e7) at (10.1,-0.4) {$E_7^{(1)}$};
\node[text=3red] (e8) at (11.9,-0.4) {$E_8^{(1)}$};

\node[eqname,text=3red,below] at (d4.south) {$PVI$};
\node[eqname,text=3red,below] at (d5.south) {$PV$};
\node[eqname,text=3red,below] at (d6.south) {$PIII^{D_6}$};
\node[eqname,text=3red,below] at (d7.south) {$PIII^{D_7}$};
\node[eqname,text=3red,below] (pd8) at (d8.south) {$PIII^{D_8}$};
\node[eqname,text=3red,below] at (e6.south) {$PIV$};
\node[eqname,text=3red,below] at (e7.south) {$PII$};
\node[eqname,text=3red,below] at (e8.south) {$PI$};

% ---------- arrows ----------
\draw[->,draw=3green] (ell) -- (m0);

\begin{scope}[->,draw=3blue]
  \draw (m0) -- (m1);  \draw (m1) -- (m2);  \draw (m2) -- (m3);
  \draw (m3) -- (m4);  \draw (m4) -- (m5);  \draw (m5) -- (m6);
  \draw (m6) -- (m7p); \draw (m6) -- (m7);  \draw (m7) -- (m8);
  \draw (m0) -- (a0);  \draw (m1) -- (a1);  \draw (m2) -- (a2);
  \draw (m3) -- (d4);  \draw (m4) -- (d5);  \draw (m5) -- (d6);
  \draw (m6) -- (d7);  \draw (m7p) -- (d8);
  \draw (m5.south west) to[out=-115,in=125] (e6.north);
  \draw (m7p.south west) -- (e7.north east);
  \draw (m8.south west) -- (e8.north east);
\end{scope}

\begin{scope}[->,draw=3ochre]
  \draw (a0) -- (a1); \draw (a1) -- (a2); \draw (a2) -- (d4);
\end{scope}

\begin{scope}[->,draw=3red]
  \draw (d4) -- (d5); \draw (d5) -- (d6); \draw (d6) -- (d7); \draw (d7) -- (d8);
  \draw (d5.south east) -- (e6.north west);
  \draw (d6.south east) -- (e7.north west);
  \draw (pd8) -- (e8);
  \draw (e6) -- (e7); \draw (e7) -- (e8);
\end{scope}
\end{tikzpicture}
\caption{Sakai's degeneration scheme. Each node is labelled by the type of the corresponding generalized Halphen surface: elliptic (green), multiplicative (blue, $q$-difference equations), additive (orange, difference equations) and by Okamoto's spaces of initial conditions of the Painlev\'e differential equations (red). Arrows indicate degenerations.}
\label{fig:sakai}
\end{figure}
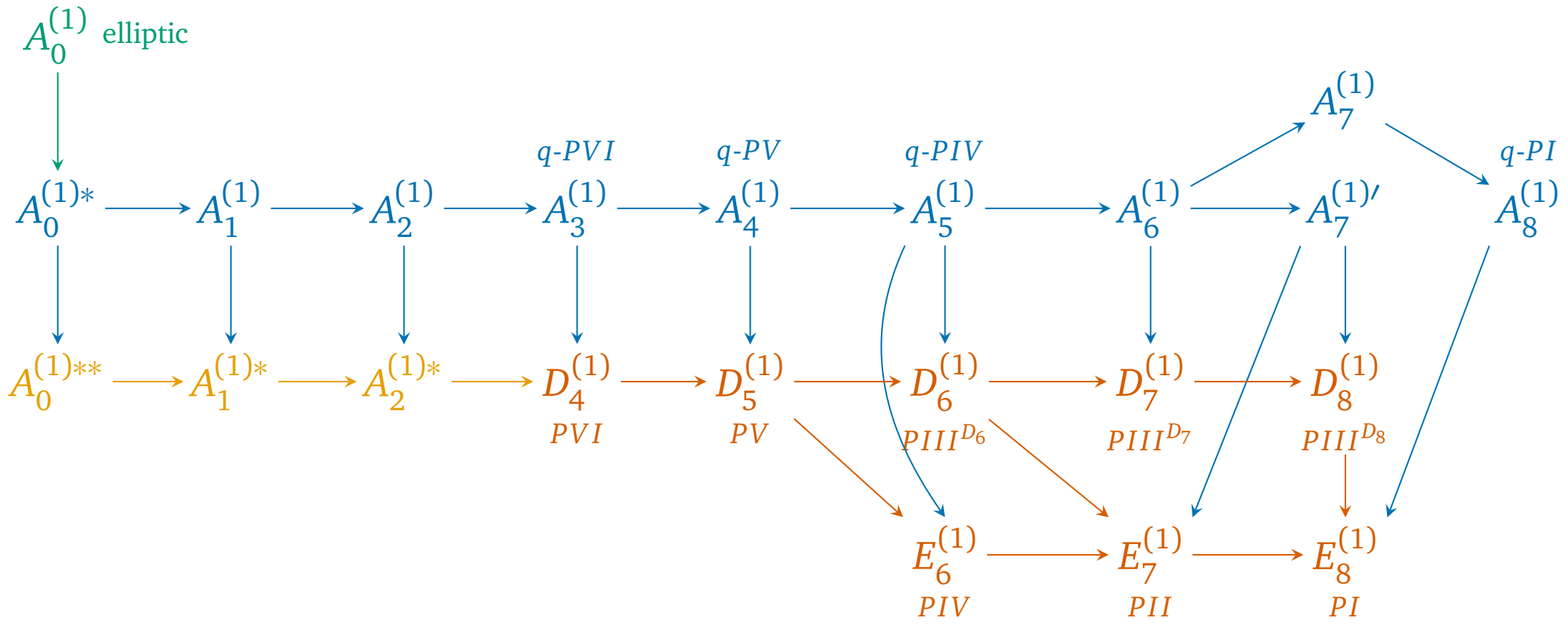

In the last decades, difference, $q$-difference and elliptic analogues of the Painlev\'e equations have been studied extensively. These are integrable discrete mappings which, in a suitable continuum limit, reduce to the Painlev\'e differential equations. In his seminal paper \cite{sakai}, Sakai gave a unified classification of the continuous and discrete Painlev\'e equations in terms of rational surfaces and their symmetries (see Figure~\ref{fig:sakai}).

Loosely speaking, each entry of Sakai's table corresponds to a Painlev\'e equation whose symmetry group is the affine Weyl group corresponding to the orthogonal lattice to the one indicated in the table (we omit some subtleties for the sake of brevity). Equations labelled in red are differential, those in orange are difference, those in blue are $q$-difference and the one in green is elliptic. The arrows represent \textit{confluences}, i.e.\ degenerations in which the independent variable, the dependent variable and the parameters are rescaled by suitable powers of a small parameter $\epsilon$ and the limit $\epsilon\to0$ is taken \cite{C}; for discrete equations, exponential changes of the dependent and independent variables may also be involved.

The Riemann--Hilbert method is the most effective approach to the global properties of the solutions of the Painlev\'e differential equations. It reconstructs solutions of the auxiliary linear system from piecewise analytic functions whose jumps across a contour $\Sigma$ are prescribed by the monodromy data. Among other results, it has been used to prove the existence of solutions that are bounded along the real line, for instance when a vanishing lemma holds for the corresponding jump conditions. For the discrete Painlev\'e equations, however, a comparable Riemann--Hilbert theory has not yet been fully developed. A natural first question in this direction is therefore to determine the monodromy manifolds of the difference, $q$-difference and elliptic Painlev\'e equations.

For the differential equations this question is well understood. By \cite{SvdP}, the monodromy manifolds $\mathcal{M}^{(d)}$ of the Painlev\'e differential equations are all affine cubic surfaces of the form
\begin{equation}\label{eq:mon-mf}
x_1 x_2 x_3 + \epsilon_1^{(d)} x_1^2+ \epsilon_2^{(d)} x_2^2+ \epsilon_3^{(d)} x_3^2 + \omega_1^{(d)} x_1 + \omega_2^{(d)} x_2 + \omega_3^{(d)} x_3+ \omega_4^{(d)}=0,
\end{equation}
where $ \epsilon_i^{(d)}\in\{0,1\}$ and
$d$ runs over $PVI$, $PV$, $PV_{\operatorname{deg}}$, $PIV$, $PIII^{D_6}$, $PIII^{D_7}$, $PIII^{D_8}$, $PII^{JM}$, $PII^{FN}$, $PI$, and different isomonodromic problems for the same equation are distinguished. The monodromy manifold of the sixth Painlev\'e equation is the $\mathrm{SL}_2(\mathbb{C})$ character variety of the four-holed Riemann sphere. Its real slice is the decorated Teichm\"uller space of the four-holed sphere, which is described combinatorially by a fat graph and shear coordinates; complexifying the shear coordinates yields flat coordinates on the full character variety \cite{ChM}. For the remaining Painlev\'e equations, the monodromy manifolds were interpreted as decorated character varieties of a Riemann sphere with boundary in \cite{CMR}. That work also gave a geometric description of the confluence procedure, the so-called \textit{chewing-gum} limit.

In view of these results, it is natural to ask whether the discrete Painlev\'e equations that degenerate to $PVI$ also have cubic monodromy manifolds of a similar form. In \cite{JR,JRM} it was shown that this expectation is too optimistic in the $q$-difference case. The present paper shows that the difference case is more favourable. In a nutshell, we show that the difference Painlev\'e equation $d$-$P\big(A_2^{(1)*}\big)$, the first discrete equation to the left of $PVI$ in Sakai's scheme, has a monodromy manifold given by a family of affine cubic surfaces, which reduces to the monodromy manifold of $PVI$ under the chewing-gum limit.

More precisely, we consider the difference Painlev\'e equation with symmetry $W\big(E_6^{(1)}\big)$, denoted $d$-$P\big(A_2^{(1)*}\big)$ in Sakai's classification. Following Dzhamay, Sakai and Takenawa \cite{DST}, this equation is realised as a composition of discrete Schlesinger transformations acting on $3\times3$ Fuchsian systems with three simple poles at $0$, $1$ and $\infty$. After a scalar normalization, the monodromy data of such a system are pairs of matrices in $\mathrm{SL}_3(\mathbb{C})$ up to simultaneous conjugation. In other words, they are points of the character variety $\mathcal{X}_{\mathrm{SL}_3}=\mathrm{SL}_3(\mathbb{C})^2/\!\!/\mathrm{SL}_3(\mathbb{C})$ of the three-holed sphere, whose coordinate ring was described by Lawton \cite{Law}.

Our first result is an explicit map from the space of local exponents of the Fuchsian system to a six-parameter family of affine cubic surfaces in $\mathbb{C}^3$ of the form
\begin{equation*}
x_1 x_2 x_3 + x_1^3+ x_2^3+ x_3^2+ \lambda_1 x_1^2+ \lambda_2 x_2^2+ \omega_1 x_1+\omega_2 x_2+\omega_3 x_3+\omega_4=0.
\end{equation*}

We show that a step of the discrete dynamics acts on $\mathcal{X}_{\mathrm{SL}_3}$ as a central twist of order three. With $\zeta=e^{2\pi i/3}$, it acts by
$$
x_1\mapsto\zeta^2x_1,\qquad x_2\mapsto\zeta x_2,\qquad x_3\mapsto x_3,
$$
and rescales the coefficients of the cubic by cube roots of unity. Every monomial of the cubic is therefore invariant, and the dynamics of the discrete Painlev\'e equation never leaves the family of cubics.

We parametrize $\mathcal{X}_{\mathrm{SL}_3}$ by a finite cover of the Fock--Goncharov torus \cite{FG,DMM}. We show that, in the chewing-gum limit in which the two coordinates attached to one edge of the triangulation tend to infinity, our cubic degenerates to the monodromy manifold of the sixth Painlev\'e equation. The traces of the four $\mathrm{SL}_2(\mathbb{C})$ monodromy matrices of $\operatorname{PVI}$ appear as explicit functions of the surviving Fock--Goncharov coordinates. This leads us to the following 

\begin{conjecture}
The monodromy manifolds of the discrete Painlev\'e equations of type $A_1^{(1)*}$ and $A_0^{(1)**}$ are respectively given by the following affine del Pezzo surfaces
$$
x_1 x_2 x_3 + x_1^4 + x_2^3 + x_3^2 + \mu_1 x_1^3 +
\lambda_1 x_1^2+ \lambda_2 x_2^2+ \omega_1 x_1+\omega_2 x_2+\omega_3 x_3+\omega_4=0
$$
$$
x_1 x_2 x_3 + x_1^5+ x_2^3+ x_3^2+ \mu_2 x_1^4 + \mu_1 x_1^3 +
\lambda_1 x_1^2+ \lambda_2 x_2^2+ \omega_1 x_1+\omega_2 x_2+\omega_3 x_3+\omega_4=0.
$$
\end{conjecture}

Another natural question for future work concerns the $q$-difference analogue of the present construction. Park \cite{Park} found a $3\times 3$ Lax pair for the $q$-Painlev\'e equation with surface type $A_2^{(1)}$. In the limit $q\to 1$, its linear problem reduces to the $3\times 3$ Fuchsian system with poles at $0,1,\infty$ studied by Dzhamay, Sakai and Takenawa \cite{DST}. This leads to the following open problem: describe the monodromy manifold of Park's $q$-difference system, and determine whether it degenerates to the monodromy manifold studied in the present paper on one side, and to the one studied by Joshi and Roffelsen \cite{JR} on the other.

This paper is organized as follows: in Section~\ref{sec2} we recall the geometry of the generalized Halphen surface of type $A_2^{(1)*}$, the standard form of the equation, and its realisation by Dzhamay, Sakai and Takenawa as discrete Schlesinger transformations of $3\times3$ Fuchsian systems. In Section~\ref{sec3} we describe the character variety $\mathcal{X}_{\mathrm{SL}_3}$ as a family of affine cubic surfaces and compute the action of the discrete dynamics on it. In Section~\ref{sec:FG} we introduce Fock--Goncharov coordinates and prove that the chewing-gum limit of our cubic is the monodromy manifold of $PVI$.

\paragraph{Acknowledgements}
The authors are grateful to 
M. Bershtein, D. Dal Martello, A. Dzhamay, I. Gaiur, N. Joshi and V. Poberezhny for helpful discussions.  The research of M. Mazzocco was supported  by the Proyecto de Generaci\'on de Conocimiento, PID2024-155686NB-I00 of the Spanish Ministry of Science, Innovation and Universities.
The research of A. Nefedova was supported by a grant of Theoretical Physics and Mathematics Advancement
Foundation “BASIS.”

\section{Fuchsian system, Schlesinger transformations, and the results of Dzhamay--Sakai--Takenawa}\label{sec2}

We first briefly recall Sakai's general theory for the equation $A_2^{(1)*}$.

\subsection{Geometry of the generalized Halphen surface of type $A_2^{(1)*}$}

In Sakai's classification \cite{sakai}, every discrete Painlev\'e equation is associated with a family of rational surfaces (namely smooth projective algebraic surfaces over $\mathbb C$ that are birationally equivalent to $\mathbb P^2$), obtained by blowing up $\mathbb{P}^2$ at nine, possibly infinitely near, points.  For the equation $d$-$P(A_2^{(1)*})$ studied here it is convenient to work instead with $\mathbb{P}^1\times\mathbb{P}^1$ blown up at eight points $p_1,\ldots,p_8$.

We start by recalling a few basic definitions and facts about this construction.

\begin{definition}\label{df:bl}
Let $X$ be a smooth algebraic surface and $p \in X$ a point. The \emph{blow-up} of $X$ at $p$ is a smooth surface $\mathrm{Bl}_p X$ together with a morphism $\pi\colon\mathrm{Bl}_p X \to X$ such that:
\begin{enumerate}
\item the \emph{exceptional divisor} $E := \pi^{-1}(p)$ is a curve isomorphic to $\mathbb{P}^1$, whose points correspond to the tangent directions to $X$ at $p$;
\item $\pi$ restricts to an isomorphism $\mathrm{Bl}_p X\setminus E \to X \setminus \{p\}$.
\end{enumerate}
In local coordinates $(x,y)$ centered at $p$, $\mathrm{Bl}_p X$ is given near $E$ by
\[
\{\, ((x,y),[u:v]) \;:\; xv = yu \,\} \subset \mathbb{A}^2 \times \mathbb{P}^1,
\]
with $\pi$ the projection to $(x,y)$ and $E = \{x=y=0\}$. The normal bundle of $E$ in $\mathrm{Bl}_p X$ is $\mathcal{O}_{\mathbb{P}^1}(-1)$, and hence $E^2 = -1$.
\end{definition}

The type of the surface is determined by how the eight points of blow-up are allowed to be positioned. Within a fixed type, the points may move, and their coordinates are described by parameters $b = \{b_1, \ldots , b_8\}$. Thus we obtain the following family of  surfaces:
$$
  X_b := \mathrm{Bl}_{p_1,\ldots,p_8}(\mathbb{P}^1\times\mathbb{P}^1).
$$

The discrete Painlevé equation defines a birational map $(f,g)\mapsto(\bar f,\bar g)$ of $\mathbb{P}^1\times\mathbb{P}^1$, together with a shift of the parameters $b\mapsto\bar b$. This map is not regular on $\mathbb{P}^1\times\mathbb{P}^1$, as it is indeterminate at the points $p_i$, but for generic parameters it lifts to an isomorphism between the corresponding blown-up surfaces:
$$
\varphi_b\colon X_b\longrightarrow X_{\bar b}.
$$

Although $\varphi_b$ is a nonlinear map, it induces a linear map $\varphi_b^*$ on the Picard group of $X_b$ that preserves the intersection form, so Sakai's classification and the identification of the dynamics as a translation in an affine Weyl group are both expressed in terms of this lattice.

\begin{definition}
The Picard group $\mathrm{Pic}(X)$ of a smooth projective surface $X$ is the group of divisors on $X$ modulo linear equivalence.
\end{definition}

For $\mathbb{P}^1\times\mathbb{P}^1$ we have
$\mathrm{Pic}(\mathbb{P}^1\times\mathbb{P}^1)
=\mathbb{Z}H_f\oplus\mathbb{Z}H_g$, where
$H_f=[f=\mathrm{const}]$ and $H_g=[g=\mathrm{const}]$.
Here $[\,\cdot\,]$ denotes the linear equivalence class: any two vertical lines $\{f=c\}$ and $\{f=c'\}$ in $\mathbb{P}^1\times\mathbb{P}^1$, viewed as divisors, satisfy
$$
\{f=c\}-\{f=c'\}=\mathrm{div}\!\left(\frac{f-c}{f-c'}\right),
$$
so they are linearly equivalent and define a single element $H_f$ of the Picard group. Similarly, all horizontal lines $\{g=c\}$ define $H_g$.
We use the same notation for the pullbacks of $H_f$ and $H_g$ to $X_b$. Each point blow-up adds a summand $\mathbb{Z}E_i$, generated by the class of the corresponding exceptional divisor, so that
$$
\mathrm{Pic}(X_b) = \mathbb{Z}H_f \oplus \mathbb{Z}H_g \oplus 
\bigoplus_{i=1}^{8} \mathbb{Z}E_i.
$$

%The Picard lattice is equipped with a symmetric bilinear form characterized by
%$$
%D \cdot C = \deg\bigl(\mathcal{O}_X(D)|_C\bigr),
%$$
%where $D$ is a divisor and $C \subset X$ is a smooth irreducible curve. In particular, if $D$ itself is a smooth irreducible curve, then
%$$
%D^2 = \deg N_{D/X}.
%$$

%For effective divisors without a common component, the intersection product is the sum of local intersection multiplicities. In our basis, the intersection pairing is given by

The Picard group carries a symmetric bilinear form, the intersection product. For two distinct irreducible curves meeting transversally it counts their intersection points, and self-intersections of exceptional curves are $E_i^2=-1$ (\Cref{df:bl}). In our basis it is given by
$$
H_f^2=H_g^2=0,\qquad H_f\cdot H_g=1,\qquad
H_f\cdot E_i=H_g\cdot E_i=0,\qquad
E_i\cdot E_j=-\delta_{ij}.
$$
Since $\mathrm{Pic}(X_b)\cong\mathbb{Z}^{10}$ is a free abelian group equipped with an integral symmetric bilinear form, the intersection product, it is common to refer to it as the \emph{Picard lattice}.

In Sakai's classification, different equations correspond to different decompositions into irreducible components, and different configurations of these components, of an effective divisor in the anti-canonical class $-K_X$, the negative of the canonical class $K_X$. The latter is the class in $\mathrm{Pic}(X)$ of the divisor of any nonzero meromorphic $2$-form $\omega$ on $X$. It does not depend on the choice of $\omega$: two such forms differ by a rational function $h$, so their divisors differ by $\mathrm{div}(h)$ and are linearly equivalent. On $\mathbb{P}^1\times\mathbb{P}^1$,
$$
\operatorname{div}(df\wedge dg)
=-2\{f=\infty\}-2\{g=\infty\}.
$$
Hence
$$
-K_{\mathbb{P}^1\times\mathbb{P}^1}=2H_f+2H_g,
\qquad
-K_{X_b}=2H_f+2H_g-\sum_{i=1}^8E_i.
$$
Moreover, $-K_{X_b}$ is preserved by the isomorphism $\varphi_b$, and it is the null root of the symmetry root lattice $R^\perp$, the orthogonal complement in $\mathrm{Pic}(X_b)$ of the root lattice $R$ that specifies the surface type of the Painlev\'e equation as in \Cref{fig:sakai}.

The surfaces $X_b$ belong to the following class of rational surfaces, which is the setting of Sakai's classification.

\begin{definition}
A smooth projective rational surface $X$ is called a generalized
Halphen surface if it admits an effective anti-canonical divisor
$$
D=\sum_i m_iD_i,\qquad [D]=-K_X,
$$
where the $D_i$ are irreducible curves, $m_i$ are positive integers,
and $K_X\cdot[D_i]=0$ for every $i$.
\end{definition}

The \emph{index} of a generalized Halphen surface $X$ is the smallest $m\geq1$ such that $\dim|-mK_X|\geq1$, or $0$ if no such $m$ exists. Surfaces of positive index carry an elliptic fibration and correspond to autonomous maps, while discrete Painlevé equations correspond to surfaces of index zero, for which $D$ is the unique effective anti-canonical divisor. 

We now describe the decomposition of the anti-canonical divisor $D$ into irreducible components and the configuration of these components for the equation $d$-$P(A_2^{(1)*})$. This equation corresponds to generalized Halphen surfaces of index zero and type $A_2^{(1)*}$. This means that the unique effective anti-canonical divisor decomposes as
$$
D = D_0 + D_1 + D_2,
\qquad [D] = -K_{X_b},
$$
where the $D_i$ are smooth rational curves with $D_i^2=-2$ and $D_i\cdot D_j=1$ for $i\neq j$.
The matrix of the intersection form on their classes is thus the negative of the affine Cartan matrix of type $A_2^{(1)}$. 
The curves $D_0,D_1,D_2$ are the strict transforms of
$f+g=0$, $f=\infty$, and $g=\infty$, respectively.
Their classes in the chosen Picard basis are
$$
[D_0]=H_f+H_g-E_1-E_2-E_3-E_4,
$$
$$
[D_1]=H_f-E_5-E_6,
$$
$$
[D_2]=H_g-E_7-E_8.
$$
The asterisk records the configuration of the components: in type $A_2^{(1)*}$ the three curves pass through a common point, $(\infty,\infty)$, while in type $A_2^{(1)}$, the
components meet pairwise at three distinct points, forming
a triangle.

The orthogonal complement of
$R=\operatorname{Span}_{\mathbb{Z}}
\{[D_0],[D_1],[D_2]\}$ in the Picard lattice is the symmetry
root lattice $R^\perp\cong Q(E_6^{(1)})$.
The corresponding affine Weyl group is $W(E_6^{(1)})$.
We choose the simple roots as follows:
$$
    \epsilon_0 = E_3-E_4,\quad
    \epsilon_1 = E_2-E_3,\quad
    \epsilon_2 = E_1-E_2,\quad
    \epsilon_3 = H_f-E_1-E_7,
$$
$$
    \epsilon_4 = E_7-E_8,\quad
    \epsilon_5 = H_g-E_1-E_5,\quad
    \epsilon_6 = E_5-E_6.
$$

\begin{center}
\begin{tikzpicture}[
    scale=0.5,
    vertex/.style={
        circle,
        draw,
        fill=white,
        minimum size=6mm,
        inner sep=0pt,
        line width=0.8pt
    },
    line width=1pt
]
    \draw (0,0) -- (8,0);
    \draw (4,0) -- (4,4);

    \node[vertex, label=below:{$\epsilon_0$}] at (0,0) {};
    \node[vertex, label=below:{$\epsilon_1$}] at (2,0) {};
    \node[vertex, label=below:{$\epsilon_2$}] at (4,0) {};
    \node[vertex, label=below:{$\epsilon_3$}] at (6,0) {};
    \node[vertex, label=below:{$\epsilon_4$}] at (8,0) {};

    \node[vertex, label=right:{$\epsilon_5$}] at (4,2) {};
    \node[vertex, label=right:{$\epsilon_6$}] at (4,4) {};
\end{tikzpicture}
\end{center}

The anti-canonical class
$$
-K_{X_b}
=\epsilon_0+2\epsilon_1+3\epsilon_2+2\epsilon_3
+\epsilon_4+2\epsilon_5+\epsilon_6
$$
is the null root of the symmetry root lattice:
it is nonzero and satisfies
$-K_{X_b}\cdot\epsilon_i=0$ for every $i$.
The induced action of the discrete evolution on the symmetry
root lattice is the translation
$$
(\epsilon_0,\epsilon_1,\epsilon_2,\epsilon_3,\epsilon_4,\epsilon_5,\epsilon_6)
\longmapsto
(\epsilon_0,\epsilon_1,\epsilon_2,\epsilon_3,\epsilon_4,\epsilon_5,\epsilon_6)
+(0,0,0,1,0,-1,0)(-K_{X_b}).
$$

\subsection{Standard form of the equation}
We use the standard form of $d$-$P(A_2^{(1)*})$ obtained by
Grammaticos, Ramani, and Ohta \cite{RGO}, with the parametrisation
of Sakai \cite{Sak07}. In coordinates $(f,g)$ on
$\mathbb{P}^1\times\mathbb{P}^1$, it takes the form

\begin{equation}\label{eq:dP}
\begin{cases}
\displaystyle
(f+g)(\bar f+g)
= \frac{(g+b_1)(g+b_2)(g+b_3)(g+b_4)}{(g-b_5)(g-b_6)},
\\
\displaystyle
(\bar f+g)(\bar f+\bar g)
= \frac{(\bar f-b_1)(\bar f-b_2)(\bar f-b_3)(\bar f-b_4)}
       {(\bar f+\bar{b}_7)(\bar f+\bar{b}_8)}.
\end{cases}
\end{equation}
One step is computed successively as
$$
(f,g)\longmapsto(\bar f,g)\longmapsto(\bar f,\bar g).
$$
The bar denotes one full step and acts on both the coordinates
and the parameters.

The parameters $b_i$ change as follows:
$$
\begin{pmatrix} 
b_1 & b_2 & b_3 & b_4 \\ 
b_5 & b_6 & b_7 & b_8 
\end{pmatrix} 
\longmapsto 
\begin{pmatrix} 
b_1 & b_2 & b_3 & b_4 \\ 
b_5 + d & b_6 + d & b_7 - d & b_8 - d
\end{pmatrix},
$$
where $d=\sum_{i=1}^8b_i$ is unchanged by the evolution.
For $d\ne0$, a simultaneous rescaling of $f,g,b_1,\ldots,b_8$
allows the normalization $d=-1$.
With the identification of the parameters in the next subsection,
$$
d=\theta_1^0+\theta_2^0+\theta_1^1+\theta_2^1
+\kappa_1+\kappa_2+\kappa_3-1=-1
$$
by the Fuchs relation. Thus this normalization imposes no
additional constraint on the exponents.

The blow-up points are
$$
\begin{aligned}
p_1&=(b_1,-b_1),& p_2&=(b_2,-b_2),\\
p_3&=(b_3,-b_3),& p_4&=(b_4,-b_4),\\
p_5&=(\infty,b_5),& p_6&=(\infty,b_6),\\
p_7&=(-b_7,\infty),& p_8&=(-b_8,\infty).
\end{aligned}
$$
The points $p_1,\ldots,p_4$ lie on $f+g=0$, whose strict
transform is $D_0$; $p_5,p_6$ lie on $f=\infty$, whose strict
transform is $D_1$; and $p_7,p_8$ lie on $g=\infty$, whose
strict transform is $D_2$.

In other words, $X_b$ is the space of initial conditions of \eqref{eq:dP}: the birational map $(f,g)\mapsto(\bar f,\bar g)$, which blows down some curves and blows up some points of $\mathbb{P}^1\times\mathbb{P}^1$, becomes bi-regular once the points $p_1,\ldots,p_8$ are blown up, and the evolution of the coordinates and the shift of the parameters are two aspects of the single isomorphism $\varphi_b$. The linear map $\varphi_b^*$ on $\mathrm{Pic}(X_b)$ essentially determines $\varphi_b$: the class $\varphi_b^*(H_f)$ prescribes the bi-degree of $\bar f$ as a rational function of $(f,g)$ and the points $p_i$ through which its level curves pass, and similarly for $\bar g$. Conversely, Sakai showed that the Cremona isometries of $\mathrm{Pic}(X_b)$, in particular the translations in $W(E_6^{(1)})$, are induced by isomorphisms between surfaces of the family \cite[Theorem 25]{sakai}. The translation in \Cref{sec2} (end of \S2.1) is the one realized by \eqref{eq:dP}.

\subsection{Fuchsian system}

The discrete Painlevé equation $d$-$P\big(A_2^{(1)*}\big)$ in the work of Dzhamay, Sakai, and Takenawa \cite{DST} is realized as dynamics on the set of Fuchsian systems of size $3\times 3$.

We consider a Fuchsian system with three simple poles at $z = 0$, $z = 1$, and $z = \infty$:
\begin{equation}\label{eq:FS3}
\frac{dY}{dz} = \left(\frac{A_0}{z} +\frac{A_1}{z - 1}\right)Y,
\end{equation}
where $A_0, A_1\in\mathrm{End}(\mathbb{C}^3)$ are diagonalizable matrices with pairwise distinct eigenvalues such that $A_\infty:= -A_0-A_1\neq 0$. By applying a scalar gauge transformation, one can assume that $A_0$ and $A_1$ have one eigenvalue equal to zero so that the Riemann scheme of any system has the form
$$ \Theta = 
\begin{pmatrix}
0 & 1 & \infty \\
\theta_1^0 & \theta_1^1 & \kappa_1 \\
\theta_2^0 & \theta_2^1 & \kappa_2 \\
0 & 0 & \kappa_3
\end{pmatrix}
$$
with the Fuchs relation $\sum_{j=1}^2 \theta_j^0 + \sum_{j=1}^2\theta_j^1 + \sum_{j=1}^3 \kappa_j = 0$.

By Katz' formula,for fixed regular semisimple residue conjugacy classes satisfying the Fuchs relation, the nonempty irreducible
moduli space of triples \((A_0,A_1,A_{\infty})\), with \(A_0+A_1+A_{\infty}=0\), modulo simultaneous conjugation has complex dimension two.
A convenient way to parameterize such a moduli space is through the decomposition of the residue matrices $A_i = B_i C_i^\dagger$; the decomposition can be chosen so that $C_i^\dagger B_i = \Theta_i = \mathrm{diag}(\theta_1^i,\theta_2^i)$, so that $B_i$ is the matrix of right eigenvectors and $C_i^\dagger$ is the matrix of left eigenvectors for the eigenvalues $\theta_1^i,\theta_2^i$. 
The third
eigenvalue of $A_i$ is $0$; the corresponding right eigenvector spans the right kernel of $C_i^\dagger$, and
the corresponding left eigenvector spans the left kernel of $B_i$.

In what follows, we denote the columns of $B_i$ by $b^i_1, b^i_2$, the columns of
$C_i$ by $c^i_1, c^i_2$, a vector spanning the right kernel of $C_i^\dagger$ by
$b^i_3$, and a vector such that $(c^i_3)^\dagger$ spans the left kernel of $B_i$ by
$c^i_3$. Thus the right eigenvectors of $A_i$ are $b^i_1, b^i_2, b^i_3$ and the left
eigenvectors are $(c^i_1)^\dagger, (c^i_2)^\dagger, (c^i_3)^\dagger$, for the
eigenvalues $\theta^i_1,\theta^i_2,0$ respectively.

By a global conjugation and rescaling, the matrices $B_i$ can be chosen as follows
$$
B_0=\begin{pmatrix}
    1 & 0 \\
    0 & 1 \\
    0 & 0
\end{pmatrix},\qquad 
B_1 =
\begin{pmatrix}
    0 & 1 \\
    0 & 1 \\
    1 & 1
\end{pmatrix}
$$
so that the matrices $C_i$ acquire the form
$$
C_0^\dagger=
\begin{pmatrix}
    \theta_1^0  & 0 & \alpha \\
    0 & \theta_2^0 & \beta 
\end{pmatrix},  \qquad
C_1^\dagger=
\begin{pmatrix}
    -\gamma - \theta_1^1  & \gamma & \theta_1^1  \\
    \theta_2^1 - \delta & \delta & 0 
\end{pmatrix},
$$
where $\alpha,\beta,\gamma,\delta$ are free parameters. 
Using the condition on the eigenvalues of the residue matrices the reduction to two coordinates  $(x,y)$ defined as 
$$
  x = \frac{(\gamma+\delta)(\theta_1^0-\theta_2^0)}{\theta_1^1-\theta_2^1},\qquad
  y = \frac{\theta_2^1\gamma + \theta_1^1\delta}{\gamma+\delta+\theta_1^1-\theta_2^1}
$$
is achieved. The formulae for $\alpha,\beta,\gamma,\delta$ in terms of $(x,y)$ can be found in \cite{DST}.

On $(x,y)$, the dynamics equivalent to the dynamics of the discrete equation is defined (see \cite[p.~24]{DST}).

\begin{remark}
    We explain here the relationship between the coordinates $(x,y)$ and the so-called \textit{spectral coordinates} $(q,p)$, namely the algebro-geometric 
Darboux coordinates constructed following the scheme of 
\cite{veselov, adams, DuDi, gekhtman}. 

Let $A(z)=\frac{A_0}{z}+\frac{A_1}{z-1}$ and denote by
\begin{equation}\label{char}
\mathcal{R}(z,w)=\det\big(w-A(z)\big)
\end{equation}
its characteristic polynomial, which has degree $3$ in $w$. The \textit{spectral curve}
\begin{equation}\label{curve}
\Gamma=\{{\cal R}(z,w)=0\}
\end{equation}
of the matrix $A(z)$ is a
$3$-sheeted cover of the $z$-sphere whose sheets are the eigenvalues
$w_1(z),w_2(z)$, $w_3(z)$ of $A(z)$. For generic values of the parameters this cover is smooth outside the
fibres over $z=0$, $z=1$ and $z=\infty$, the poles of $A(z)$. Over each
of these points the (compactified) curve $\Gamma$ has exactly three
points, and the residues of the differential $w\,dz$ at them are
$\{\theta^0_1,\theta^0_2,0\}$, $\{\theta^1_1,\theta^1_2,0\}$ and
$\{\kappa_1,\kappa_2,\kappa_3\}$ respectively, i.e.\ they are given by
the Riemann scheme $\Theta$.

In order to fix the matrix $A(z)$ uniquely (up to global conjugation), we need to select one more point $(q,p)\in\Gamma$. The procedure consists in selecting this point in a smart way as explained in \cite{skly}.

Let us introduce the row vectors ${v}_0$, ${v}_1(z)$ and $v_2(z)$,  by
\begin{equation}\label{bk}
v_0=(0, 0,1), ~~v_1(z) = v_0 A(z),  ~~v_2(z) = v_1(z) A(z)
\end{equation}
Denote $V(z)$ the $3\times 3$ matrix with the rows $v_0$, $v_1(z)$, 
$v_{2}(z)$. Put
\begin{equation}\label{Delta}
{\Delta_0}(z) =  z^3(z-1)^3  
\det V(z).
\end{equation}
It is easy to see that the equation ${\Delta_0}(z)=0$ has two roots, $z=0$ (double) and $z=q$ where $q:=\frac{x}{x+y}$. At $z=q$ the vectors $v_0,v_1(q),v_2(q)$ are therefore linearly
dependent and, for generic $(x,y)$, there exists a unique, up to normalization, nonzero vector $\psi$ spanning $\ker V(q)$. Then:
\begin{itemize}
\item $v_0\psi=0$, i.e.\ the third component of $\psi$ vanishes;
\item $\psi$ is an eigenvector of
$A(q)$, and we call $p$ the corresponding eigenvalue:
\[
A(q)\psi=p\,\psi,\qquad p=\frac{(x+y)(\theta^0_2-x)}{x}.
\]
\end{itemize}
In particular $\mathcal R(q,p)=0$, i.e.\ $(q,p)$ is a point of the
spectral curve. Conversely, the coordinates $(x,y)$ are recovered from
$(q,p)$ by
\[
x=\theta^0_2-pq,\qquad y=\frac{(1-q)(\theta^0_2-pq)}{q}.
\]
\end{remark}

\subsection{Dynamics on the set of Fuchsian systems}
The equation $d$-$P\big(A_2^{(1)*}\big)$ is obtained by successively applying two transformations to a Fuchsian system. One of them is an elementary Schlesinger transformation.

\begin{definition}
An elementary Schlesinger transformation {\footnotesize{$\left\{ \begin{array}{l}\tau \phi \\ \mu \nu \end{array} \right\}$} }
is a gauge transformation $\tilde{Y} (z) = R(z)Y(z)$, where the matrix $R(z)$ has the form
$$R(z) = I + \frac{z_{\tau} - z_{\phi}}{z - z_{\tau}} P, $$
$P$ is a rank-one projector determined by eigenvectors $\mu, \nu$ of the residue matrices. Under such a transformation, the poles do not change, but the eigenvalue $\theta_{\mu}^{\tau}$ decreases by 1, while $\theta_{\nu}^{\phi}$ increases by 1; all other eigenvalues remain unchanged.
\end{definition}

Since $R(z)$ is rational in $z$, it is single-valued and therefore preserves the monodromy representation exactly. In addition, one uses the map $\Sigma_0(1,3) = \rho_0(\theta_1^0)\circ\sigma_0(1,3)$, which consists of two steps:

\begin{definition}
The transformation $\rho_i(s)$ acts on the residue matrices by
$$
A_i \longmapsto A_i - sI, \qquad A_\infty \longmapsto A_\infty + sI,
$$
leaving all other matrices unchanged. At the level of the Riemann scheme, it shifts all exponents at the point $z_i$ by $-s$ and all exponents at $\infty$ by $+s$.
\end{definition}

The transformation $\rho_i(s)$ is induced by the scalar gauge transformation
$\tilde Y(z)=(z-z_i)^{-s}Y(z)$. Hence it multiplies the monodromy matrices by scalars,
$$
\tilde M_i = e^{-2\pi\mathrm{i}s}M_i,\qquad \tilde M_\infty = e^{2\pi\mathrm{i}s}M_\infty,
$$
and leaves the other monodromy matrices unchanged. In particular, it preserves the monodromy exactly when $s\in\mathbb{Z}$.

\begin{definition}
The transformation $\sigma_i(j,k)$ exchanges the $j$-th eigenpair
$(b^i_j,c^i_j)$ of the matrix $A_i$ with the $k$-th eigenpair
$(b^i_k,c^i_k)$, together with the corresponding exponents
$\theta^i_j\leftrightarrow\theta^i_k$. 
\end{definition}

Note that the transformation $\sigma_i(j,k)$  does
 not change any of the matrices $A_0,A_1,A_\infty$. It changes only the labeling of the eigenvalues and eigenvectors of $A_i$. Therefore the monodromy matrices are also preserved under this transformation.

The action of $\Sigma_0(1,3)$ at the level of the Riemann scheme looks as follows:
$$
  (\theta_1^0,\theta_2^0,0)
  \xrightarrow{\sigma_0(1,3)}
  (0,\theta_2^0,\theta_1^0)
  \xrightarrow{\rho_0(\theta_1^0)}
  (-\theta_1^0,\theta_2^0-\theta_1^0,0),
  \qquad \kappa_j\longmapsto\kappa_j+\theta_1^0.
$$

The parameters $b_i$ determining the Okamoto surface are related to the parameters of the Fuchsian system as follows:
$$
b_1 = \theta^0_2+\kappa_1,\quad
b_2 = \theta^0_2+\kappa_2,\quad
b_3 = \theta^0_2+\kappa_3,\quad
b_4 = 0,
$$
$$
b_5 = \theta^1_1,\quad
b_6 = \theta^1_2,\quad
b_7 = \theta^0_1-\theta^0_2,\quad
b_8 = -\theta^0_2-1.
$$
\begin{theorem}[Dzhamay, Sakai, Takenawa \cite{DST}]\label{thm:dP}
The equation $d$-$P\big(A_2^{(1)*}\big)$ is given as the composition of the following transformations:
$$ \Sigma_{0}(1,3)\circ 
\left\{ \begin{array}{l}1 0\\ 2 1 \end{array} \right\} \circ \Sigma_{0}(1,3)\circ \left\{ \begin{array}{l}1 0\\ 1 1 \end{array} \right\}.$$
This composition acts on the Riemann scheme as
$$(\theta_1^0,\theta_2^0,\theta_1^1,\theta_2^1,\kappa_1,\kappa_2,\kappa_3)\longmapsto (\theta_1^0,\theta_2^0 -1,\theta_1^1 -1,\theta_2^1 -1,\kappa_1 + 1,\kappa_2 + 1,\kappa_3 + 1).$$
This is an integer shift $\Theta \longmapsto \Theta +\delta$, where $\delta = (0, - 1, - 1, - 1, + 1, + 1, + 1)$.
\end{theorem}

\section{The $\mathrm{SL}_3(\mathbb{C})$ character variety}\label{sec3}

Section~\ref{sec2} describes the relation between the dynamics on the set of
Fuchsian systems with $3\times 3$ matrices and the dynamics of the discrete Painlevé
equation. Therefore we can study the corresponding monodromy-data variety and verify
that the discrete dynamics on the set of Fuchsian systems preserves this variety.

The monodromy of system~\eqref{eq:FS3} defines a representation
$$
\rho:\pi_1\big(\overline{\mathbb{C}}\setminus\{0,1,\infty\}\big)\longrightarrow
\mathrm{GL}_3(\mathbb{C}).
$$
Since the fundamental group of the sphere with three punctures is a free group on two
generators, such a representation is determined by the monodromy matrices $M_0,M_1$,
and $M_\infty=(M_0M_1)^{-1}$. Hence
$\mathrm{Hom}\big(\pi_1(\overline{\mathbb{C}}\setminus\{0,1,\infty\}),G\big)\cong G^2$
for any group $G$, and monodromy data up to simultaneous conjugation are described by
the following variety.

\begin{definition}
Let $G$ be a complex reductive algebraic group. The $G$-character variety of the
sphere with three holes is the affine algebraic variety
$$
\mathcal{X}_G = G^2 /\!\!/ G,
$$
with coordinate ring
$$
\mathbb{C}[\mathcal{X}_G] = \mathbb{C}[G\times G]^{G},
$$
where $G$ acts on pairs $(A,B)\in G^2$ by simultaneous conjugation,
$g\cdot(A,B)=(gAg^{-1},gBg^{-1})$.
\end{definition}

The monodromy matrices of~\eqref{eq:FS3} lie in $\mathrm{GL}_3(\mathbb{C})$, and the
$\mathrm{GL}_3(\mathbb{C})$-character variety is related to the
$\mathrm{SL}_3(\mathbb{C})$-character variety by the isomorphism
$$
\big(\mathcal{X}_{\mathrm{SL}_3}\times(\mathbb{C}^*)^2\big)\big/\mathbb{Z}_3^2
\;\xrightarrow{\ \sim\ }\;
\mathcal{X}_{\mathrm{GL}_3}
$$
induced by the map
$\big([A,B],(\lambda_0,\lambda_1)\big)\mapsto[\lambda_0 A,\lambda_1B]$. Here the group
$\mathbb{Z}_3^2$ acts on $\mathcal{X}_{\mathrm{SL}_3}$ by central twists,
$$
(k_0,k_1)\cdot[A,B]=[\zeta^{k_0}A,\zeta^{k_1}B],\qquad \zeta=e^{2\pi\mathrm{i}/3},
$$
and on $(\mathbb{C}^*)^2$ by
$(k_0,k_1)\cdot(\lambda_0,\lambda_1)=(\zeta^{-k_0}\lambda_0,\zeta^{-k_1}\lambda_1)$.
In particular, setting $\widehat M_i=e^{-2\pi\mathrm{i}\operatorname{tr}(A_i)/3}M_i$,
the monodromy of~\eqref{eq:FS3} corresponds to the pair
$(A,B)=(\widehat M_0,\widehat M_1)\in\mathrm{SL}_3(\mathbb{C})^2$, with
$\widehat M_\infty=(AB)^{-1}$.

\begin{theorem}
The character variety $\mathcal{X}_{\mathrm{SL}_3} = \mathrm{SL}_3(\mathbb{C})^2 
\!\mathbin{/\!/} \mathrm{SL}_3(\mathbb{C})$ is isomorphic to an affine hypersurface of degree $6$ in $\mathbb{C}^9$.
The coordinate ring $\mathbb{C}[\mathcal{X}_{\mathrm{SL}_3}]$ is generated by the following trace functions:
\begin{equation}\label{eq:tracesAB}
\begin{aligned}
t(1) &= \operatorname{tr}(A), \qquad & t(-1) &= \operatorname{tr}(A^{-1}),\\
t(2) &= \operatorname{tr}(B), \qquad & t(-2) &= \operatorname{tr}(B^{-1}),\\
t(3) &= \operatorname{tr}(AB), \qquad & t(-3) &= \operatorname{tr}((AB)^{-1}),\\
t(4) &= \operatorname{tr}(AB^{-1}),\qquad & t(-4) &= \operatorname{tr}(A^{-1}B),\\
t(5) &= \operatorname{tr}([A,B]) = \operatorname{tr}(ABA^{-1}B^{-1}), 
\end{aligned}
\end{equation}
subject to the cubic relation 
\begin{equation}\label{eq:cubic}
x_1 x_2 x_3 + x_1^3+ x_2^3+ x_3^2+ \lambda_1 x_1^2+ \lambda_2 x_2^2+
\omega_1 x_1+\omega_2 x_2+\omega_3 x_3+\omega_4=0,
\end{equation}
where the variables $x_1,x_2,x_3$ are defined as
\begin{equation}\label{eq:shift}
\begin{split}
&
x_{1} =  t(4) - t(1)t(-2),\qquad
x_{2} =  t(-4) - t(-1)t(2),\\
&
x_{3} = -t(5)+t(1)t(-1) + t(2)t(-2) + t(3)t(-3) - 6,
\end{split}
\end{equation}
and the parameters are expressed through the traces $t(\pm 1),t(\pm 2), t(\pm 3)$ that fix the conjugacy classes of the monodromy matrices, as follows:
\begin{equation} \label{eq:lambdas_omega3}
\begin{aligned}
\lambda_1 &= t(1)t(-2) + t(-1)t(-3) + t(2)t(3),
\\
\lambda_2 &= t(1)t(3) + t(-1)t(2) + t(-2)t(-3),
\\
\omega_3 &= - t(-1) t(-2) t(3) - t(-1) t(1) - t(-2) t(2) - t(-3) t(1) t(2) - t(-3) t(3) + 9,
\end{aligned}
\end{equation}
\begin{equation} \label{eq:omega1}
\begin{aligned}
\omega_1 = & t(-1)^2 t(3) + t(-1) t(-2)^2 + t(-1) t(-2) t(-3) t(1) + t(-1) t(-3) t(2) t(3) -\\
& - 3 t(-1) t(2) - 3 t(-2) t(-3) + t(-2) t(1) t(2) t(3) + t(-2) t(3)^2 +\\
& + t(-3)^2 t(1) + t(-3) t(2)^2 + t(1)^2 t(2) - 3 t(1) t(3),\\
\end{aligned}
\end{equation}
\begin{equation} \label{eq:omega2}
\begin{aligned}
\omega_2 = & t(-1)^2 t(-2) + t(-1) t(-2) t(-3) t(2) - 3 t(-1) t(-3) + t(-1) t(1) t(2) t(3) +\\
& + t(-1) t(3)^2 + t(-2)^2 t(3) + t(-2) t(-3) t(1) t(3) - 3 t(-2) t(1) +\\
& + t(-3)^2 t(2) + t(-3) t(1)^2 + t(1) t(2)^2 - 3 t(2) t(3),\\
\end{aligned}
\end{equation}

\begin{equation} \label{eq:omega4}
\begin{aligned}
\omega_4 = &t(-1)^3 + t(-1)^2 t(-2)^2 t(-3) + t(-1)^2 t(-2) t(1) t(3) - 2 t(-1)^2 t(-3) t(2) +\\
& + t(-1)^2 t(2) t(3)^2 + t(-1) t(-2)^2 t(2) t(3) - 2 t(-1) t(-2) t(-3)^2 + t(-1) t(-2) t(-3) t(1) t(2) t(3) +\\
& + t(-1) t(-2) t(-3) t(3)^2 + t(-1) t(-2) t(1) t(2) - 3 t(-1) t(-2) t(3) + t(-1) t(-3)^2 t(2)^2 +\\
&+ t(-1) t(-3) t(1)^2 t(2) + t(-1) t(-3) t(1) t(3) - 9 t(-1) t(1) - 2 t(-1) t(2)^2 t(3)+\\
&+ t(-2)^3 - 2 t(-2)^2 t(-3) t(1) + t(-2)^2 t(1) t(3)^2 + t(-2) t(-3)^2 t(1)^2+\\
& +t(-2) t(-3) t(1) t(2)^2 + t(-2) t(-3) t(2) t(3) - 2 t(-2) t(1)^2 t(3) - 9 t(-2) t(2) +\\
&+ t(-3)^3 + t(-3)^2 t(1) t(2) t(3) - 3 t(-3) t(1) t(2) - 9 t(-3) t(3)+\\
&+ t(1)^3 + t(1)^2 t(2)^2 t(3) - 2 t(1) t(2) t(3)^2 + t(2)^3+\\
&+ t(3)^3 + 27.\\
\end{aligned}
\end{equation}
The singular locus of $\mathcal{X}_{\mathrm{SL}_3}$ coincides with the image of reducible representations of pairs of matrices having a common eigenvector.
\end{theorem}

\begin{proof}
The fact that the character variety $\mathcal{X}_{\mathrm{SL}_3}$ is isomorphic to an affine hypersurface of degree $6$ in $\mathbb{C}^9$ and that its singular locus coincides with the image of reducible representations
is due to Lawton \cite{Law}. He also proved that  the coordinate ring $\mathbb{C}[\mathcal{X}_{\mathrm{SL}_3}]$ is generated by the trace functions \eqref{eq:tracesAB} subject to 
the equation
\begin{equation}\label{eq:Lawton}
t(5)^{2} - P t(5) + Q = 0,
\end{equation}
where 
\begin{equation}\label{eq:P}
\begin{aligned}
P &= t(1)t(-1)t(2)t(-2) - t(1)t(2)t(-3) - t(-1)t(-2)t(3)\\
&\quad - t(1)t(-2)t(-4) - t(-1)t(2)t(4) \\
&\quad + t(1)t(-1) + t(2)t(-2) + t(3)t(-3) + t(4)t(-4) - 3,
\end{aligned}
\end{equation}
\begingroup
\allowdisplaybreaks
\begin{align}
Q &= 9 - 6t(1)t(-1) - 6t(2)t(-2) - 6t(3)t(-3) - 6t(4)t(-4)
    + t(1)^3 + t(2)^3 + t(3)^3 \notag\\
&\quad + t(4)^3 + t(-1)^3 + t(-2)^3 + t(-3)^3 + t(-4)^3
    - 3t(-4)t(-3)t(-1) - 3t(4)t(3)t(1) \notag\\
&\quad - 3t(-4)t(2)t(3) - 3t(4)t(-2)t(-3)
    + 3t(-4)t(-2)t(1) + 3t(4)t(2)t(-1) \notag\\
&\quad + 3t(1)t(2)t(-3) + 3t(-1)t(-2)t(3)
    + t(-2)t(-1)t(2)t(1) + t(-3)t(-2)t(3)t(2) \notag\\
&\quad + t(-4)t(-1)t(4)t(1) + t(-4)t(-2)t(4)t(2)
    + t(-3)t(-1)t(3)t(1) + t(-3)t(-4)t(3)t(4) \notag\\
&\quad + t(-4)^2t(-3)t(-2) + t(4)^2t(3)t(2)
    + t(-1)^2t(-2)t(-4) + t(1)^2t(2)t(4) + t(1)t(-2)^2t(-3) \notag\\
&\quad + t(-1)t(2)^2t(3) + t(-4)t(-3)t(1)^2
    + t(4)t(3)t(-1)^2 + t(-4)t(2)t(-3)^2 + t(4)t(-2)t(3)^2 \notag\\
&\quad + t(-1)^2t(-3)t(2) + t(1)^2t(3)t(-2)
    + t(-4)t(1)t(2)^2 + t(4)t(-1)t(-2)^2 + t(-4)t(3)t(-2)^2 \notag\\
&\quad + t(4)t(-3)t(2)^2 + t(1)t(3)t(-4)^2
    + t(-1)t(-3)t(4)^2 + t(-1)t(-4)t(3)^2 \notag\\
&\quad + t(1)t(4)t(-3)^2 - 2t(-3)^2t(-2)t(-1)
    - 2t(3)^2t(2)t(1) - 2t(-4)^2t(-1)t(2) \notag\\
&\quad - 2t(4)^2t(1)t(-2) + t(-1)^2t(-2)^2t(-3)
    + t(1)^2t(2)^2t(3) + t(-4)t(-1)^2t(2)^2 \notag\\
&\quad + t(4)t(1)^2t(-2)^2 - t(-4)t(-2)^2t(2)t(1)
    - t(4)t(2)^2t(-2)t(-1) \notag\\
&\quad - t(-3)t(1)^2t(-1)t(2) - t(3)t(-1)^2t(1)t(-2)
    - t(-3)t(2)^2t(-2)t(1) - t(3)t(-2)^2t(2)t(-1) \notag\\
&\quad - t(-4)t(-2)t(-1)t(1)^2 - t(4)t(2)t(1)t(-1)^2
    - t(-1)t(-2)^3t(1) - t(-1)t(2)^3t(1) \notag\\
&\quad - t(-1)^3t(-2)t(2) - t(1)^3t(-2)t(2)
    - t(-4)t(-3)t(-2)t(-1)t(2) - t(4)t(3)t(2)t(1)t(-2) \notag\\
&\quad - t(-1)t(1)t(2)t(-4)t(3) - t(-1)t(1)t(-2)t(4)t(-3)
    + t(-2)^2t(-1)^2t(2)t(1) \notag\\
&\quad + t(-1)^2t(-2)t(2)^2t(1). \label{eq:Q}
\end{align}
\endgroup
In equation \eqref{eq:Lawton}, the six traces 
 $t(\pm1), t(\pm2), t(\pm3)$ fix the conjugacy classes of the monodromy matrices, while the remaining three variables $t(\pm4), t(5)$ parameterize the internal degrees of freedom. The change of coordinates 
\eqref{eq:shift} was found by looking for linear transformations for $t(\pm 4)$ that map the part depending on $t(\pm 4)$ in \eqref{eq:Lawton} 
$$
t(4)t(-4) - t(-1)t(2)t(4) - t(1)t(-2)t(-4),
$$
to the correct form in $x_1,x_2$. Then the transformation was completed for $x_3$.  
\end{proof}

The discrete dynamics of
Theorem~\ref{thm:dP} preserves the monodromy matrices $M_0,M_1,M_\infty\in\mathrm{GL}_3(\mathbb{C})$  up to simultaneous conjugation, hence every point of $\mathcal{X}_{\mathrm{GL}_3}$
is fixed by the dynamics. At the same time, the traces of the residue matrices change
by
$$
(\operatorname{tr}A_0,\operatorname{tr}A_1,\operatorname{tr}A_\infty)\longmapsto
(\operatorname{tr}A_0-1,\ \operatorname{tr}A_1-2,\ \operatorname{tr}A_\infty+3).
$$
Because the traces shift while the $M_i$ are
fixed, the normalized monodromy matrices
$$
\widehat M_i=e^{-2\pi\mathrm{i}\operatorname{tr}(A_i)/3}M_i
$$
transform as
$$
(\widehat M_0,\widehat M_1,\widehat M_\infty)\longmapsto
(\zeta\widehat M_0,\ \zeta^2\widehat M_1,\ \widehat M_\infty),
\qquad \zeta=e^{2\pi\mathrm{i}/3}.
$$
Thus, with $A=\widehat M_0$ and $B=\widehat M_1$, the discrete dynamics acts on
$\mathcal{X}_{\mathrm{SL}_3}$ by the central twist
$$
[A,B]\longmapsto[\zeta A,\ \zeta^2 B],
$$
which leaves $\widehat M_\infty=(AB)^{-1}$ unchanged. Since $\zeta\,\mathbb{I}$ is central
in $\mathrm{SL}_3(\mathbb{C})$, this map commutes with simultaneous conjugation and
defines an automorphism of $\mathcal{X}_{\mathrm{SL}_3}$ of order $3$. In terms of the
isomorphism $\big(\mathcal{X}_{\mathrm{SL}_3}\times(\mathbb{C}^*)^2\big)/\mathbb{Z}_3^2
\cong\mathcal{X}_{\mathrm{GL}_3}$, with $\lambda_i=e^{2\pi\mathrm{i}\operatorname{tr}(A_i)/3}$,
the dynamics acts on $\mathcal{X}_{\mathrm{SL}_3}\times(\mathbb{C}^*)^2$ by the element
$(1,2)\in\mathbb{Z}_3^2$, and therefore trivially on $\mathcal{X}_{\mathrm{GL}_3}$. 

On the trace
coordinates, the discrete dynamic acts by
\begin{equation}\label{eq:tr-twist}
\begin{aligned}
\operatorname{tr}A &\mapsto \zeta\operatorname{tr}A, &\qquad
\operatorname{tr}A^{-1} &\mapsto \zeta^2\operatorname{tr}A^{-1},\\
\operatorname{tr}B &\mapsto \zeta^2\operatorname{tr}B, &\qquad
\operatorname{tr}B^{-1} &\mapsto \zeta\operatorname{tr}B^{-1},\\
\operatorname{tr}(AB^{-1}) &\mapsto \zeta^2\operatorname{tr}(AB^{-1}), &\qquad
\operatorname{tr}(A^{-1}B) &\mapsto \zeta\operatorname{tr}(A^{-1}B),
\end{aligned}
\end{equation}
while $\operatorname{tr}(AB)$, $\operatorname{tr}\big((AB)^{-1}\big)$ and
$\operatorname{tr}[A,B]$ are invariant. In particular, the dynamics preserves
$\mathcal{X}_{\mathrm{GL}_3}$, and every orbit in $\mathcal{X}_{\mathrm{SL}_3}$ consists of at most three points.

Substituting the transformations \eqref{eq:tr-twist} into formulas \eqref{eq:lambdas_omega3}-\eqref{eq:omega4} and \eqref{eq:shift}, we obtain the following transformations of the parameters of the cubic:
\begin{align*}
&
\lambda_1 \longmapsto e^{-2\pi i / 3} \lambda_1, \quad
\lambda_2 \longmapsto e^{2\pi i / 3} \lambda_2, \\
&\omega_1 \longmapsto e^{2\pi i / 3} \omega_1, \quad
\omega_2 \longmapsto e^{-2\pi i / 3} \omega_2, \quad
\omega_3 \longmapsto \omega_3, \quad
\omega_4 \longmapsto  \omega_4, \\
\end{align*}
and on the coordinates of the cubic by
\begin{equation}\label{eq:x-twist}
x_1\longmapsto \zeta^{2}x_1=e^{-2\pi i/3}x_1,\qquad
x_2\longmapsto \zeta\, x_2=e^{2\pi i/3}x_2,\qquad
x_3\longmapsto x_3,
\end{equation}
so that, together with the transformation of $\lambda_1,\lambda_2,\omega_1,\dots,\omega_4$ given above, every monomial in \eqref{eq:cubic} is invariant and the cubic is preserved by the dynamics.

\section{Fock--Goncharov coordinates and the chewing-gum limit}\label{sec:FG}

In this section we parametrize the character variety $\mathcal X_{\mathrm{SL}_3}$ of Section \ref{sec3}  by a finite cover of the Fock–Goncharov torus \cite{FG} and use it to produce a limit in which the cubic \eqref{eq:cubic} degenerates to the cubic \eqref{eq:mon-mf} of the sixth Painlev\'e equation.

\subsection{Fock--Goncharov parametrization}

\begin{figure}[h]
    \centering
    \includegraphics[width=.6\textwidth]{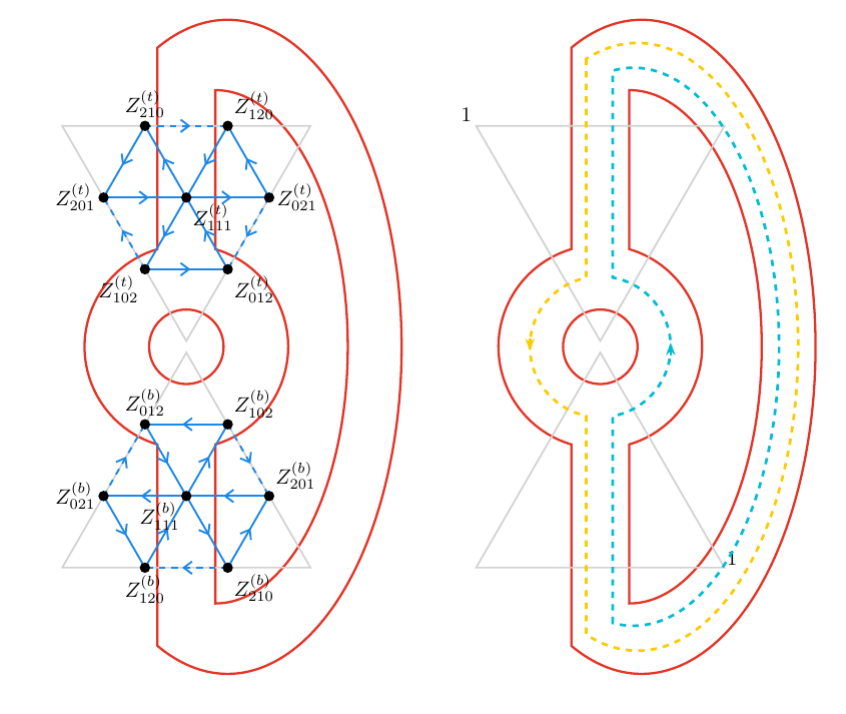}
    \caption{The ribbon graph of the three-holed Riemann sphere in red. On the left, the dual triangulation is displayed in gray with the corresponding Fock--Goncharov variables. On the right the paths corresponding to the monodromy matrices $A$ and $B$ are shown in yellow and cyan respectively.}
    \label{fig:triangulation}
\end{figure}

The parametrization of the monodromy data of the Fuchsian system was given in  \cite[Section~4]{DMM}. We recall parts of that construction here. 
The three-holed sphere can be described by a planar ribbon graph with three holes (see red graph in \Cref{fig:triangulation}, in which two of the holes are inside and the third is outside). It can also be triangulated by two  triangles, which are dual to the ribbon graph, in \Cref{fig:triangulation} they are given by a top one and a bottom one, glued along their three edges. According to Fock-Goncharov prescription, in the rank $3$ case, each triangle carries seven coordinates $Z_{ijk}$, $i+j+k=3$: two on each side, namely $(Z_{201},Z_{102})$, $(Z_{120},Z_{210})$, $(Z_{012},Z_{021})$, and one, $Z_{111}$, in the interior. We denote by $Z^{t}_{ijk}$ and $Z^{b}_{ijk}$ the variables of the top and of the bottom triangle respectively. 

The monodromy matrices $A$ and $B$ correspond to loops in the ribbon graph and can be factorized in terms of a matrix $S$ and transition matrices (see right hand side of \Cref{fig:triangulation}). Denote
$$
S=\begin{pmatrix} 0&0&1\\ 0&-1&0\\ 1&0&0\end{pmatrix},
$$
$$
T_1=\left(\frac{Z_{201}Z_{102}^2}{Z_{120}Z_{210}^2Z_{111}}\right)^{\frac13}
\begin{pmatrix}
1 & Z_{210}\left(1+Z_{111}\right) & Z_{111}Z_{210}Z_{120}\\[2pt]
-Z_{102}^{-1} & -Z_{210}Z_{102}^{-1} & 0\\[2pt]
Z_{102}^{-1}Z_{201}^{-1} & 0 & 0
\end{pmatrix},
$$
and
$$
T_2=\left(\frac{Z_{120}Z_{210}^2}{Z_{012}Z_{021}^2Z_{111}}\right)^{\frac13}
\begin{pmatrix}
1 & Z_{021}\left(1+Z_{111}\right) & Z_{111}Z_{021}Z_{012}\\[2pt]
-Z_{210}^{-1} & -Z_{021}Z_{210}^{-1} & 0\\[2pt]
Z_{210}^{-1}Z_{120}^{-1} & 0 & 0
\end{pmatrix}.
$$
Gluing the two triangles along the sides that are connected by edges in the ribbon graph amalgamates the side variables in pairs as follows:
\begin{equation}\label{eq:amalg}
\begin{aligned}
Z_{Y_1}&=Z^{t}_{201}Z^{b}_{021},&\qquad Z_{Y_2}&=Z^{t}_{102}Z^{b}_{012},&\qquad Z_{Y_3}&=Z^{t}_{210}Z^{b}_{120},\\
Z_{C_1}&=Z^{t}_{021}Z^{b}_{201},&\qquad Z_{C_2}&=Z^{t}_{012}Z^{b}_{102},&\qquad Z_{C_3}&=Z^{t}_{120}Z^{b}_{210}.
\end{aligned}
\end{equation}
Together with the two face variables $Z^{t}_{111}$ and $Z^{b}_{111}$ this gives eight coordinates, which is the dimension of $\mathcal X_{\mathrm{SL}_3}$.

The monodromy matrices are then conjugated to the matrices corresponding to the 
 paths in \Cref{fig:triangulation}:
\begin{equation}\label{eq:FG-AB}
A = C\, S\, T_2^{\,b}\, S\, T_1^{\,t}\, C^{-1},\qquad
B = C\, T_2^{\,t}\, S\, T_1^{\,b}\, S\, C^{-1},
\end{equation}
where $T_i^{\,t}$ and $T_i^{\,b}$ denote $T_i$ evaluated on the variables of the top and of the bottom triangle, and
$$
C=\diag\left(\big(Z^{t}_{210}\big)^{-\frac23}\big(Z^{t}_{120}\big)^{-\frac13},\
\big(Z^{t}_{210}\big)^{\frac13}\big(Z^{t}_{120}\big)^{-\frac13},\
\big(Z^{t}_{210}\big)^{\frac13}\big(Z^{t}_{120}\big)^{\frac23}\right)
$$
is a diagonal conjugation chosen so that, after imposing \eqref{eq:amalg}, the entries of $A$ and $B$ are expressed through the amalgamated variables.
Explicitly:
\begin{multline*}
A=\frac{1}{\left(Z_{Y_1}^2Z_{Y_2}Z_{Y_3}^2Z_{C_3}Z^{t}_{111}Z^{b}_{111}\right)^{\frac13}}
\left(\begin{matrix}
1\vphantom{\Bigl(Z^{b}_{111}}\\[4pt]
Z_{Y_3}\left(1+Z_{Y_1}\right)\vphantom{\Bigl(Z^{b}_{111}}\\[4pt]
Z_{Y_3}Z_{C_3}\Bigl(1+Z_{Y_1}\bigl(1+Z^{b}_{111}(1+Z_{Y_2})\bigr)\Bigr)
\end{matrix}\right.\\
\left.\begin{matrix}
0 & 0\vphantom{\Bigl(Z^{b}_{111}}\\[4pt]
Z_{Y_1}Z_{Y_3} & 0\vphantom{\Bigl(Z^{b}_{111}}\\[4pt]
Z_{Y_1}Z_{Y_3}Z_{C_3}\Bigl(1+Z^{b}_{111}\bigl(1+Z_{Y_2}(1+Z^{t}_{111})\bigr)\Bigr)
& Z_{Y_1}Z_{Y_2}Z_{Y_3}Z_{C_3}Z^{t}_{111}Z^{b}_{111}
\end{matrix}\right),
\end{multline*}

\begin{multline*}
B=\frac{1}{\left(Z_{C_1}^2Z_{C_2}Z_{C_3}^2Z_{Y_3}Z^{t}_{111}Z^{b}_{111}\right)^{\frac13}}
\left(\begin{matrix}
Z_{C_1}Z_{C_2}Z_{C_3}Z_{Y_3}Z^{t}_{111}Z^{b}_{111}\vphantom{\Bigl(Z^{b}_{111}}\\[4pt]
0\vphantom{\Bigl(Z^{b}_{111}}\\[4pt]
0\vphantom{\Bigl(Z^{b}_{111}}
\end{matrix}\right.\\
\left.\begin{matrix}
-Z_{C_1}Z_{C_3}\Bigl(1+Z^{t}_{111}\bigl(1+Z_{C_2}(1+Z^{b}_{111})\bigr)\Bigr)
& 1+Z_{C_1}\bigl(1+Z^{t}_{111}(1+Z_{C_2})\bigr)\\[4pt]
Z_{C_1}Z_{C_3} & -1-Z_{C_1}\vphantom{\Bigl(Z^{b}_{111}}\\[4pt]
0 & 1\vphantom{\Bigl(Z^{b}_{111}}
\end{matrix}\right).
\end{multline*}
In particular the eigenvalues of $A$, $B$ and $AB$ are monomials in the Fock--Goncharov variables:
\begin{align*}
\operatorname{Spec}(A)&=\Bigl\{
Z_{Y_1}^{-\frac23}Z_{Y_2}^{-\frac13}Z_{Y_3}^{-\frac23}Z_{C_3}^{-\frac13}Z_{111}^{t^{-\frac13}}Z_{111}^{b^{-\frac13}},\ 
Z_{Y_1}^{\frac13}Z_{Y_2}^{-\frac13}Z_{Y_3}^{\frac13}Z_{C_3}^{-\frac13}Z_{111}^{t^{-\frac13}}Z_{111}^{b^{-\frac13}},\ 
Z_{Y_1}^{\frac13}Z_{Y_2}^{\frac23}Z_{Y_3}^{\frac13}Z_{C_3}^{\frac23}Z_{111}^{t^{\frac23}}Z_{111}^{b^{\frac23}}
\Bigr\},\\
\operatorname{Spec}(B)&=\Bigl\{
Z_{C_1}^{-\frac23}Z_{C_2}^{-\frac13}Z_{C_3}^{-\frac23}Z_{Y_3}^{-\frac13}Z_{111}^{t^{-\frac13}}Z_{111}^{b^{-\frac13}},\ 
Z_{C_1}^{\frac13}Z_{C_2}^{-\frac13}Z_{C_3}^{\frac13}Z_{Y_3}^{-\frac13}Z_{111}^{t^{-\frac13}}Z_{111}^{b^{-\frac13}},\ 
Z_{C_1}^{\frac13}Z_{C_2}^{\frac23}Z_{C_3}^{\frac13}Z_{Y_3}^{\frac23}Z_{111}^{t^{\frac23}}Z_{111}^{b^{\frac23}}
\Bigr\},\\
\operatorname{Spec}(AB)&=\Bigl\{
Z_{Y_1}^{-\frac23}Z_{Y_2}^{-\frac13}Z_{C_1}^{-\frac23}Z_{C_2}^{-\frac13}Z_{111}^{t^{-\frac23}}Z_{111}^{b^{-\frac23}},\ 
Z_{Y_1}^{\frac13}Z_{Y_2}^{-\frac13}Z_{C_1}^{\frac13}Z_{C_2}^{-\frac13}Z_{111}^{t^\frac13}Z_{111}^{b^\frac13},\ 
Z_{Y_1}^{\frac13}Z_{Y_2}^{\frac23}Z_{C_1}^{\frac13}Z_{C_2}^{\frac23}Z_{111}^{t^\frac13}Z_{111}^{b^\frac13}
\Bigr\}.
\end{align*}

Note that the prefactors of $A$ and $B$, which implement the normalization
from $\mathrm{GL}_3(\mathbb{C})$ to $\mathrm{SL}_3(\mathbb{C})$ of
Section~\ref{sec3}, involve cubic roots of the Fock--Goncharov variables.
Hence $A$ and $B$ are functions on the $\mathbb{Z}_3^2$-cover of the Fock--Goncharov torus
$(\mathbb{C}^*)^8$. The coordinates and the coefficients of the cubic
\eqref{eq:cubic} are invariant under the diagonal twist
$[A,B]\mapsto[\zeta A,\zeta B]$, and therefore only involve the ratio of
the two prefactors,
\begin{equation}\label{eq:varrho}
\varrho=\left(\frac{Z_{Y_1}^2Z_{Y_2}Z_{Y_3}}
{Z_{C_1}^2Z_{C_2}Z_{C_3}}\right)^{\frac13}.
\end{equation}
Thus the family of cubics \eqref{eq:cubic} is parametrised by a triple
cover of the Fock--Goncharov torus, and the discrete dynamics
\eqref{eq:x-twist} is its deck transformation
$\varrho\mapsto\zeta\varrho$.

\begin{proposition}\label{prop:FG}
Let $A,B$ be given by \eqref{eq:FG-AB}, then all parameters defined in \eqref{eq:lambdas_omega3}, \eqref{eq:omega1}, \eqref{eq:omega2} and \eqref{eq:omega4},  as well as the variables 
$x_1,x_2,-x_3$ defined in \eqref{eq:shift} are Laurent polynomials with positive coefficients in terms of the Fock--Goncharov variables 
$Z_{C_1},Z_{C_2},Z_{C_3},Z_{Y_1},Z_{Y_2},Z_{Y_3},Z^{t}_{111},Z^{b}_{111}$ and the cube root $\varrho$ defined in \eqref{eq:varrho}. Moreover, the $\lambda_1,\lambda_2,\omega_1,\omega_2,\omega_3,\omega_4,x_1,x_2,x_3$ parameterized in this way 
satisfy the cubic relation \eqref{eq:cubic} identically in $Z_{Y_1},Z_{Y_2},Z_{Y_3},Z_{C_1},Z_{C_2},Z_{C_3},Z^{t}_{111},Z^{b}_{111}$.
\end{proposition}

\begin{proof}
     The formulae are rather large and can be found in the Mathematica file \cite{MF}, where the cubic relation is also verified.
\end{proof}

\subsection{The chewing-gum limit}

We now perform the chewing-gum operation on the rightmost edge in \Cref{fig:triangulation}, namely, we set
\begin{equation}\label{eq:gum}
Z_{Y_3}=\varepsilon^{-1}\widetilde Z_{Y_3},\qquad Z_{C_3}=\varepsilon^{-1}\widetilde Z_{C_3},
\end{equation}
and let $\varepsilon\to0$ keeping $Z_{Y_1},Z_{Y_2},Z_{C_1},Z_{C_2},Z^{t}_{111},Z^{b}_{111},\widetilde Z_{Y_3},\widetilde Z_{C_3}$ fixed.

\begin{proposition}\label{prop:gum}
Under \eqref{eq:gum} the following limits exist:
$$
\xi_1=\lim_{\varepsilon\to0}\varepsilon\,x_1,\qquad
\xi_2=\lim_{\varepsilon\to0}\varepsilon\,x_2,\qquad
\xi_3=\lim_{\varepsilon\to0}\varepsilon^2x_3,
$$
$$
\nu_i=\lim_{\varepsilon\to0}\varepsilon^2\lambda_i,\quad i=1,2,\qquad
\varpi_j=\lim_{\varepsilon\to0}\varepsilon^3\omega_j,\quad j=1,2,\qquad
\varpi_3=\lim_{\varepsilon\to0}\varepsilon^2\omega_3,\qquad
\varpi_4=\lim_{\varepsilon\to0}\varepsilon^4\omega_4,
$$
and they satisfy
\begin{equation}\label{eq:lim-cubic}
\xi_1\xi_2\xi_3+\xi_3^2+\nu_1\xi_1^2+\nu_2\xi_2^2
+\varpi_1\xi_1+\varpi_2\xi_2+\varpi_3\xi_3+\varpi_4=0.
\end{equation}
Explicitly,
$$
\xi_1=\left(\frac{Z_{C_2}Z_{Y_1}\widetilde Z_{Y_3}^{\,2}\widetilde Z_{C_3}}{Z_{C_1}Z_{Y_2}}\right)^{\frac13}P_1,\qquad
\xi_2=\left(\frac{Z_{C_1}Z_{Y_2}\widetilde Z_{Y_3}\widetilde Z_{C_3}^{\,2}}{Z_{C_2}Z_{Y_1}}\right)^{\frac13}P_2,\qquad
\xi_3=-\widetilde Z_{Y_3}\widetilde Z_{C_3}Z^{t}_{111}P_3,
$$
where
\begin{align*}
P_1&=1+Z^{b}_{111}\left(1+Z_{C_1}\right)\left(1+Z_{Y_2}\right)+Z_{C_1}Z_{Y_2}Z^{t}_{111}Z^{b}_{111},\\
P_2&=1+Z^{t}_{111}\left(1+Z_{C_2}\right)\left(1+Z_{Y_1}\right)+Z_{C_2}Z_{Y_1}Z^{t}_{111}Z^{b}_{111},\\
P_3&=\left(1+Z_{C_2}\right)\left(1+Z_{Y_1}\right)
+Z^{b}_{111}\Big(Z_{C_2}+Z_{Y_1}+Z_{C_2}Z_{Y_1}\left(2+Z_{C_1}+Z_{Y_2}\right)\Big)\\
&\quad+\big(Z^{b}_{111}\big)^2Z_{C_2}Z_{Y_1}\left(1+Z_{C_1}\right)\left(1+Z_{Y_2}\right).
\end{align*}
and $\nu_1,\nu_2$ are the monomials
$$
\nu_1=\left(Z_{C_1}^2Z_{C_2}Z_{Y_1}Z_{Y_2}^2\widetilde Z_{Y_3}^{\,2}\widetilde Z_{C_3}^{\,4}\right)^{\frac13}Z^{t}_{111}Z^{b}_{111},\qquad
\nu_2=\left(Z_{C_1}Z_{C_2}^2Z_{Y_1}^2Z_{Y_2}\widetilde Z_{Y_3}^{\,4}\widetilde Z_{C_3}^{\,2}\right)^{\frac13}Z^{t}_{111}Z^{b}_{111}.
$$
\end{proposition}

\begin{proof}
Multiplying \eqref{eq:cubic} by $\varepsilon^4$, every term has a finite limit except $\varepsilon^4x_1^3=\varepsilon\,(\varepsilon x_1)^3$ and $\varepsilon^4x_2^3=\varepsilon\,(\varepsilon x_2)^3$, which vanish as $\varepsilon\to0$. The explicit expressions follow by direct computation and can be found in the mathematica file \cite{MF}.
\end{proof}

In order to match the cubic \eqref{eq:lim-cubic} with the one associated to the sixth Painlev\'e equation, which is of the form \eqref{eq:mon-mf} with all $\epsilon_i^{(d)}=1$, we need to perform the following rescaling 
$$
(\nu_1,\nu_2,\varpi_1,\varpi_2,\varpi_3,\varpi_4)\longmapsto
\left(\frac{\nu_1}{b^2},\frac{\nu_2}{a^2},\frac{\varpi_1}{ab^2},\frac{\varpi_2}{a^2b},\frac{\varpi_3}{ab},\frac{\varpi_4}{a^2b^2}\right),
$$
with $a^2=\nu_2$ and $b^2=\nu_1$.
This is equivalent to fixing the finite parts in \eqref{eq:gum} as
\begin{equation}\label{eq:norm}
\widetilde Z_{Y_3}=\left(Z_{C_2}Z_{Y_1}Z^{t}_{111}Z^{b}_{111}\right)^{-\frac12},\qquad
\widetilde Z_{C_3}=\left(Z_{C_1}Z_{Y_2}Z^{t}_{111}Z^{b}_{111}\right)^{-\frac12}.
\end{equation}
In conclusion, we have proved the following:
\begin{theorem}\label{thm:gum}
With the normalization \eqref{eq:norm}, the chewing-gum limit of the cubic \eqref{eq:cubic} is
\begin{equation}\label{eq:PVI-cubic}
\xi_1\xi_2\xi_3+\xi_1^2+\xi_2^2+\xi_3^2
+\varpi_1\xi_1+\varpi_2\xi_2+\varpi_3\xi_3+\varpi_4=0,
\end{equation}
where
\begin{equation}\label{eq:xi}
\xi_1=\frac{P_1}{\left(Z_{C_1}Z_{Y_2}Z^{t}_{111}Z^{b}_{111}\right)^{\frac12}},\qquad
\xi_2=\frac{P_2}{\left(Z_{C_2}Z_{Y_1}Z^{t}_{111}Z^{b}_{111}\right)^{\frac12}},\qquad
\xi_3=-\frac{P_3}{\left(Z_{C_1}Z_{C_2}Z_{Y_1}Z_{Y_2}\right)^{\frac12}Z^{b}_{111}},
\end{equation}
and
\begin{equation}\label{eq:varpi}
\begin{aligned}
\varpi_1&=G_1G_\infty+G_2G_3,\qquad
\varpi_2=G_2G_\infty+G_1G_3,\qquad
\varpi_3=-G_3G_\infty-G_1G_2,\\
\varpi_4&=G_1^2+G_2^2+G_3^2+G_\infty^2+G_1G_2G_3G_\infty-4,
\end{aligned}
\end{equation}
with
\begin{equation}\label{eq:G}
\begin{aligned}
G_1&=\left(Z_{C_2}Z_{Y_2}Z^{t}_{111}Z^{b}_{111}\right)^{\frac12}+\left(Z_{C_2}Z_{Y_2}Z^{t}_{111}Z^{b}_{111}\right)^{-\frac12},&\qquad
G_3&=\left(Z_{Y_1}Z_{Y_2}\right)^{\frac12}+\left(Z_{Y_1}Z_{Y_2}\right)^{-\frac12},\\
G_2&=\left(Z_{C_1}Z_{Y_1}Z^{t}_{111}Z^{b}_{111}\right)^{\frac12}+\left(Z_{C_1}Z_{Y_1}Z^{t}_{111}Z^{b}_{111}\right)^{-\frac12},&\qquad
G_\infty&=\left(Z_{C_1}Z_{C_2}\right)^{\frac12}+\left(Z_{C_1}Z_{C_2}\right)^{-\frac12}.
\end{aligned}
\end{equation}
\end{theorem}

After the change of sign $(\xi_1,\xi_2,\xi_3)\to(-\xi_1,-\xi_2,\xi_3)$, equation \eqref{eq:PVI-cubic} with parameters \eqref{eq:varpi} is the monodromy manifold of the sixth Painlev\'e equation, i.e. the $\mathrm{SL}_2(\mathbb C)$ character variety of the four-holed Riemann sphere, $G_1,G_2,G_3,G_\infty$ being the traces of the monodromy matrices around the four holes \cite{ChM,CMR}. The six variables $Z_{Y_1},Z_{Y_2},Z_{C_1},Z_{C_2},Z^{t}_{111},Z^{b}_{111}$ surviving the limit thus give the four parameters \eqref{eq:G} of the cubic and a parametrization \eqref{eq:xi} of the two-dimensional surface \eqref{eq:PVI-cubic}.

\begin{remark}
    The action of the discrete dynamics does not survive the chewing-gum limit, because it is absorbed by the normalization \eqref{eq:norm}. Indeed, by \eqref{eq:tr-twist} the limits of Proposition~\ref{prop:gum} transform as \begin{equation}
\begin{aligned}
(\xi_1,\xi_2,\xi_3) &\mapsto (\zeta^2\xi_1,\zeta\xi_2,\xi_3),\\
(\nu_1,\nu_2) &\mapsto (\zeta^2\nu_1,\zeta\nu_2),\\
(\varpi_1,\varpi_2,\varpi_3,\varpi_4) &\mapsto (\zeta\varpi_1,\zeta^2\varpi_2,\varpi_3,\varpi_4),
\end{aligned}
\end{equation}
which is precisely the rescaling above with $(a,b)=(\zeta,\zeta^2)$. Hence, once $\nu_1$ and $\nu_2$ are normalised to $1$, the cubic \eqref{eq:PVI-cubic} and its parameters \eqref{eq:varpi} are invariant under the dynamics, up to the choice of branch of the square roots in \eqref{eq:norm}.
\end{remark}

We conclude by explaining the match between the limiting cubic
\eqref{eq:PVI-cubic} and the monodromy manifold of $PVI$ in terms of
\cite{DMM}: the chewing-gum limit \eqref{eq:gum} removes the two vertices
$Z_{Y_3},Z_{C_3}$ from the diamond-shaped quiver of the rank $3$
Fock--Goncharov variables of the three-holed sphere, leaving the same full
subquiver on $Z_{Y_1},Z_{Y_2},Z_{C_1},Z_{C_2},Z^{t}_{111},Z^{b}_{111}$ that
is produced by the two quiver seizures of \cite[Section~5]{DMM}, and by
\cite[Theorem~16]{DMM} a cluster mutation at $Z^{b}_{111}$ followed by a
monomial map identifies this subquiver with the quiver of the four-holed
sphere, the explicit change of variables to the exponentiated shear
coordinates $s_1,s_2,s_3$ and perimeters $p_1,p_2,p_3$ of \cite{ChM,CMR}
being
$$
e^{s_1}=Z_{Y_1}\big(1+Z^{b}_{111}\big),\quad
e^{s_2}=Z^{t}_{111},\quad
e^{s_3}=\frac{Z_{C_1}Z^{b}_{111}}{1+Z^{b}_{111}},\quad
e^{p_1}=\frac{Z_{C_2}}{Z_{Y_1}},\quad
e^{p_2}=\frac{1}{Z^{t}_{111}Z^{b}_{111}},\quad
e^{p_3}=\frac{Z_{Y_2}}{Z_{C_1}},
$$
under which the monomials appearing in \eqref{eq:G} are Casimir functions.

\end{document}